# COMPUTER-INTENSIVE RATE ESTIMATION, DIVERGING STATISTICS AND SCANNING

By Tucker McElroy and Dimitris N. Politis

*U.S. Bureau of the Census and University of California, San Diego*

A general rate estimation method is proposed that is based on studying the in-sample evolution of appropriately chosen diverging/converging statistics. The proposed rate estimators are based on simple least squares arguments, and are shown to be accurate in a very general setting without requiring the choice of a tuning parameter. The notion of scanning is introduced with the purpose of extracting useful subsamples of the data series; the proposed rate estimation method is applied to different scans, and the resulting estimators are then combined to improve accuracy. Applications to heavy tail index estimation as well as to the problem of estimating the long memory parameter are discussed; a small simulation study complements our theoretical results.

**1. Introduction.** Let $X_1, \ldots, X_n$ be an observed stretch from a general time series $\{X_t\}$ that is not necessarily linear, or stationary. A number of converging and/or diverging statistics can be computed from a dataset of this type. In many instances, however, the rate of convergence/divergence of some statistics of interest may be unknown, that is, it may depend on some unknown feature of the underlying probability law $P$. This rate is often a quantity of direct interest; for example, it may be connected to the heavy tail index, the long memory or self-similarity parameter, and so on.

For each given context, that is, choice of statistic and assumptions on the time series $\{X_t\}$, a context-specific rate estimator may be devised and its properties analyzed. By contrast, a general approach for rate estimation has been given in the subsampling literature where knowledge/estimation of the rate is necessary for the construction of confidence intervals, hypothesis tests, and so on; see Bertail, Politis and Romano [3] or Politis, Romano and Wolf [19], Chapter 8. The subsampling rate estimator is based on evaluating the statistic of interest over subsamples of different size; subsequently, the rate of convergence/divergence is gauged by the effect incurred on the distribution of the statistic when the subsample size varies.

The subsampling rate estimator is consistent under very weak conditions. Nevertheless, a typical condition assumed in connection with subsampling is the strong mixing condition which may preclude its applicability in settings where the time series exhibits long-range dependence. In addition, the subsample size must be carefully chosen for optimal results; in general, this is a difficult problem analogous to the notorious bandwidth choice problem in nonparametric smoothing; see Politis, Romano and Wolf [19], Chapter 9.

In this paper, a different noncontext-specific rate estimation method is introduced based on studying the in-sample evolution of appropriately chosen converging/diverging statistics. The proposed rate estimator is based on a simple least squares argument and is shown to be consistent in a very general setting that does not require the strong mixing assumption. Furthermore, no "bandwidth-type" selection is required for the new estimator.

In order to improve the accuracy of this general estimation method, the notion of scanning a sequence is introduced. The proposed rate estimation method is implemented over different "scans" of the data sequence $X_1, \ldots, X_n$, and the resulting estimators are then combined to yield an improved estimator in the spirit of the "bagging" aggregation of Breiman [4].

In the next section a motivating example is given in the setup of estimation of the heavy tail index with data from a linear time series model. Section 3 introduces the general rate estimation methodology based on statistics that converge/diverge without centering; the important notion of scanning a sequence is also introduced. In Section 4 the methodology is extended to cover the case of statistics that require centering in order to converge. An application to the problem of estimating the long memory parameter of a long-range dependent time series is given in Section 5, together with a novel application combining heavy tails and long-range dependence. The setup of Section 2 is revisited in Section 6 by means of a finite-sample simulation; all proofs are deferred to the Appendix.

## 2. A motivating example: the heavy tail index.

2.1. *A heavy-tailed linear time series.* Throughout this section (and this section only) we will assume that the data $X_1, \ldots, X_n$ are an observed stretch of a *linear* time series satisfying $X_t = \sum_{j \in \mathbf{Z}} \psi_j Z_{t-j}$, for all $t \in \mathbf{Z}$, where $\{Z_t\}$ is i.i.d. from some distribution $F \in D(\alpha)$. The filter coefficients $\{\psi_j\}$ are assumed to be absolutely summable, and $D(\alpha)$ denotes the domain of attraction of an $\alpha$-stable law with $\alpha \in (0, 2]$; see, for example, Embrechts, Klüppelberg and Mikosch [8], Chapter 2.

In this context, it is well known that there exist sequences $a_n$ and $b_n$ such that $a_n^{-1}(\sum_{t=1}^n Z_t - b_n) \stackrel{\mathcal{L}}{\Longrightarrow} S_\alpha$, where $S_\alpha$ denotes a generic $\alpha$-stable law with unspecified scale, location and skewness; recall that $a_n = n^{1/\alpha} \tilde{L}(n)$ for some slowly-varying function $\tilde{L}(\cdot)$. The centering sequence $b_n$ can be taken to be zero if either $\alpha < 1$ or $\alpha > 1$ and $Z_t$ has mean zero. When $\alpha = 1$, we can only let $b_n = 0$ if $Z_t$ is symmetric about zero.

Our goal is estimation of $\alpha$, which is tantamount to estimation of the main part of the rate $a_n$; the shape of the unknown slowly-varying function $\tilde{L}(\cdot)$ is thus considered a nuisance parameter. Tail index estimators typically are based upon a number $q$ of extreme order statistics, such as the well-known Hill estimator; see Csörgő, Deheuvels and Mason [5] and Csörgő and Viharos [6]. A practical problem for these estimators is choosing the number of order statistics, such $q$ to be used; while it is known that we must have $q \to \infty$ and $q/n \to 0$ as $n \to \infty$ to ensure consistency, the optimal choice of $q$ in any given finite sample situation is challenging; see, for example, Danielsson et al. [7] and the references therein.

An alternative tail index estimator that is not based on order statistics has been recently proposed in the subsampling literature; see Bertail, Politis and Romano [3] and Politis, Romano and Wolf [19], Chapter 8. The subsampling tail index estimator is consistent under very general conditions; interestingly, it shares with Hill's estimator the difficulty of having to choose a "bandwidth"-type parameter, namely, the subsample size. It is of interest to construct a general rate estimator that is free from this difficulty of a "bandwidth"-type selection.

2.2. *A simple tail index estimator.* Let $S_n^2 = \frac{1}{n} \sum_{t=1}^n X_t^2$, and note that when $\alpha \in (0, 2)$ it follows that

$$(1) \qquad n^{-2/\alpha} L(n) \sum_{t=1}^n X_t^2 \xrightarrow{\mathcal{L}} J,$$

where $L(\cdot)$ is a slowly-varying function and $J$ has a positively skewed $S_{\alpha/2}$ distribution; see, for example, McElroy and Politis [14]. When $\alpha = 2$, the expression (1) is valid if $Z_t$ has finite variance, and the law of large numbers kicks in.

Let $Y_k = \log S_k^2$, and $U_k = Y_k - \gamma \log k + \log L(k)$ for $k = 1, \ldots, n$, where $\gamma = -1 + 2/\alpha$. Then it is immediate that (1) implies that $U_n = O_P(1)$. From the relation $Y_k = \gamma \log k + U_k - \log L(k)$, it is suggested that $\gamma$ could plausibly be estimated as the slope of a regression of $Y_k$ on $\log k$, with a resulting estimator for $\alpha$. The reason that we treat $\log L(k)$ as "approximately constant" in the regression of $Y_k$ on $\log k$ is given in Proposition 2.1 below.

PROPOSITION 2.1. *Any slowly-varying function $L(\cdot)$ satisfies*

$$(2) \qquad \log L(k) = o(\log k) \qquad as \ k \to \infty.$$

So, let $\hat{\gamma}$ and $\check{\gamma}$ be the slope estimators in least squares (LS) regression of $Y_k$ on $\log k$ with and without an intercept term, respectively, and define $\check{\alpha} = 2/(\check{\gamma} + 1)$ and $\hat{\alpha} = 2/(\hat{\gamma} + 1)$. A rough estimate of slope in the regression without an intercept is simply the ratio $Y_n/\log n$; see Meerschaert and Scheffler [16].

PROPOSITION 2.2. *If* (1) *is true, then* $\check{\alpha} \xrightarrow{P} \alpha$, *as* $n \to \infty$.

Proposition 2.2—whose proof follows from the more general Theorem 3.1 in the next section—remains true even in the case $\alpha = 2$ as long as (1) holds; see the discussion in Section 6.

The rate of convergence of $\check{\alpha}$ can be quite slow. To get a more accurate estimator, a permutation/averaging technique was proposed in Politis [17]. However, permutations are only justified in the special case when the $X_t$ data are i.i.d.; to address the general scenario of dependent data, the notion of *scanning* is introduced in Section 3.2 and will be used in connection with an estimator of the type of $\hat{\alpha}$. Intuitively, including an intercept term in the regression offers an improvement, as it captures the nonzero large-sample expectation of $U_k$, as well as the influence of the term $\log L(k)$.

## 3. The general rate estimation methodology.

3.1. *Statistics that converge or diverge without centering.* We outline below the basic rate estimation method and show its consistency under general conditions.

(a) Let $T_n = T_n(X_1, \ldots, X_n)$ be some positive statistic whose rate of convergence/divergence depends on some unknown real-valued parameter $\lambda$.

(b) Assume that for some slowly varying function $L(n)$ and for some known invertible function $g(\cdot)$ that is continuous over an interval that contains $\lambda$, we have $U_n = O_P(1)$ as $n \to \infty$, where

$$(3) \qquad U_k = \log(k^{-g(\lambda)} L(k) T_k) \qquad \text{for } k = 1, \ldots, n.$$

(c) Estimate $g(\lambda)$ by $\hat{g} = \frac{\sum_{k=1}^n (Y_k - \bar{Y})(\log k - \overline{\log n})}{\sum_{k=1}^n (\log k - \overline{\log n})^2}$, and $\lambda$ by $\hat{\lambda} = g^{-1}(\hat{g})$, where $Y_k = \log T_k$ for $k = 1, \ldots, n$, $\bar{Y} = \frac{1}{n} \sum_{k=1}^n Y_k$ and $\overline{\log n} = \frac{1}{n} \sum_{k=1}^n \log k$. Alternatively, estimate $g(\lambda)$ by $\check{g} = \sum_{k=1}^n Y_k \log k / \sum_{k=1}^n \log^2 k$, and $\lambda$ by $\check{\lambda} = g^{-1}(\check{g})$.

To study $\hat{\lambda}$ and $\check{\lambda}$, the following additional assumptions will be useful:

$$(4) \qquad U_n \xRightarrow{\mathcal{L}} \text{some r.v. } U, \qquad \text{with } EU_n^2 \to EU^2,$$

$$(5) \qquad EU_n - EU = O(n^{-p}) \qquad \text{for some } p > 0,$$

RATE ESTIMATION AND SCANNING 5and

(6) $\quad \text{Cov}(U_b, U_n) = O(b^{\gamma_1} n^{-\gamma_2} \tilde{L}(n))$ for $b \leq n$ and some $0 \leq \gamma_1 < \gamma_2$,

where $\tilde{L}$ is some slowly-varying function. Equations (4), (5) and/or (6) can be verified under some assumptions in the setting of Section 2; see www.math.ucsd.edu/~politis/PAPER/scansAppendix2.pdf for details.

We are now able to state a general asymptotic result on $\check{\lambda}$ and $\hat{\lambda}$. Theorem 3.1 below is a more general (and corrected) version of results in Politis [17] that were worked out under the assumption that the slowly-varying function is a constant.

THEOREM 3.1. *Assume statements* (a), (b), (c) *are true.*

(i) *Then* $\check{\lambda} \xrightarrow{P} \lambda$ *as* $n \to \infty$.
(ii) *If assumption (4) holds, then* $E\hat{g} \to g(\lambda)$ *and* $\text{Var}(\hat{g}) = O(1)$ *as* $n \to \infty$.
(iii) *If assumption (5) holds, then* $E\hat{g} = g(\lambda) + A_1 + A_2$, *where*

$$A_1 = O(n^{-p} \log n) \quad \text{and} \quad A_2 = -\frac{\sum_{k=1}^n (\log L(k) - \overline{\log L})(\log k - \overline{\log n})}{\sum_{k=1}^n (\log k - \overline{\log n})^2}.$$

(iv) *If assumptions (4) and (6) hold, then* $\text{Var}(\hat{g}) = o(1)$ *and* $\hat{\lambda} \xrightarrow{P} \lambda$.

REMARK 3.1. Note that the estimators $\hat{g}$ and $\check{g}$ correspond to $L_2$ regression estimators of slope (with or without an intercept). However, an $L_1$ regression estimator of slope would be a robust alternative which is expected to also be consistent and perhaps even more reliable, especially if the large-sample distribution of the $U_k$ has heavier tails than the normal.

REMARK 3.2. The assumption $U_n = O_P(1)$ in statement (b) would typically be verified by proving a limit theorem of the type

(7) $\qquad n^{-g(\lambda)} L(n) T_n \xRightarrow{\mathcal{L}} J \qquad \text{as } n \to \infty,$

where $J$ is some well-defined probability distribution. Therefore, the implication of the assumption $U_n = O_P(1)$ is that if $g(\lambda) > 0$, then $T_n$ diverges to $\infty$, whereas if $g(\lambda) < 0$, then $T_n$ converges to 0 in probability; the case $g(\lambda) = 0$ roughly corresponds to the case where the uncentered distribution of $T_n$ converges in law to some nondegenerate distribution. Unless $g(\lambda) = 0$, $Y_k = \log T_k$ diverges to either $+\infty$ or $-\infty$ as the block size $k$ increases. In addition, note that centering can typically be omitted only when $T_n$ is a *diverging* statistic, in particular, when the centering is constant or grows at a slower rate than the scale of $T_n$. Thus, most applications of Theorem 3.1 are expected to be in cases where $g(\lambda) \geq 0$. However, this rule is not adamant, as Remark 3.3 suggests.

REMARK 3.3. In the setting of Section 2 the parameter $\lambda$ would be the heavy tail index $\alpha$, and $T_n$ could well be the second sample moment $S_n^2$; in that case, $g(\lambda) = -1 + 2/\lambda$. Note, however, that the diverging statistic $S_n^2$ can be turned into a statistic that converges to zero by appropriate shrinking. For example, if $\alpha > 1$, then the statistic $T_n' = S_n^2/n$ converges weakly to zero, and $\log T_n'$ to $-\infty$; thus, the choice of $U_n$ based on the statistic $T_n'$ ensuring $U_n = O_P(1)$ is *identical* to the $U_n$ corresponding to the diverging statistic $S_n^2$, and $\hat{\alpha}$ is the *same* in both cases, which is reassuring. In essence, these are not really separate cases; since $\log(n^d T_n) = d \log n + \log T_n$, multiplying the statistic $T_n$ by $n^d$ leads to the same log–log regression.

REMARK 3.4. The validity of the regression of $Y_k$ on $\log k$ is based on asymptotic assumptions such as $U_n = O_P(1)$ or $U_n \stackrel{\mathcal{L}}{\Longrightarrow} U$, line (2), and so on. Hence, the $(Y_k, \log k)$ points may not be very informative if $k$ is small, and it may be advisable in practice to drop some points from the regression, much in the same manner as some points are invariably dropped in the beginning of a Markov chain simulation. In other words, one would regress $Y_k$ on $\log k$ for $k = n_0, \ldots, n$, for some $n_0$ chosen either as constant or even as a function of $n$ but such that $n - n_0 \to \infty$ without affecting the asymptotic consistency of $\breve{\lambda}$ or $\hat{\lambda}$. Thus, choosing $n_0$ here is not a bandwidth-choice problem, and the choice $n_0 = 1$ is definitely a valid one; the reason is that the log–log scatterplot is very sparse for points with $k$ small, and therefore, such points have little influence collectively.

Theorem 3.1 shows that $\breve{\lambda}$ is consistent under minimal assumptions, essentially the $U_n = O_P(1)$ assumption of statement (b). Nevertheless, the rate of convergence of $\breve{\lambda}$ may be very slow, essentially of logarithmic order. Intuitively, as mentioned in Section 2, the estimator $\hat{\lambda}$ should be more accurate than $\breve{\lambda}$; this is indeed true at the expense of the additional assumptions (4), (5) and (6). For example, it is immediate that the bias of $\hat{\lambda}$ will tend to zero at a polynomial rate under some conditions on the slowly-varying function $L$, for example, when $L$ is constant. However, no rate for the variance of $\hat{\lambda}$ was given in Theorem 3.1. Furthermore, if assumption (6) fails and/or can not be verified, the rough bound $\mathrm{Var}(\hat{\lambda}) = O(1)$ ensues by the delta method. Therefore, a technique to reduce the variance of $\hat{\lambda}$ is desirable; this is accomplished in the next subsection via the notion of scanning a sequence.

3.2. *Scanning a sequence.* The rate estimation method introduced in Section 3.1 is based on evaluating the statistic $T_k$ on subsets (blocks) of growing size taken from the data set $X_1, \ldots, X_n$. Subsequently, the in-sample evolution of the (logarithm of the) statistic $T_k$ is studied. This method is closely related to subsampling since our statistic is evaluated on subsamples/subseries of the data. The only difference is that here we consider blocks of all sizes as opposed to one preferred block size; as a matter of fact, here we have one block for each block size $k = 1, \ldots, n$. As in subsampling, the crux of the method outlined in Section 3.1 lies in the fact that $T_k$ and $T_{k'}$ should behave similarly (when appropriately normalized); see Politis and Romano [18] or Politis, Romano and Wolf [19] for more details on the subsampling methodology, and Barbe and Bertail [1] in connection with the study of subsamples of increasing size.

To fix ideas, assume that the time series $\{X_t\}$ is strictly stationary. In that case, it is apparent that the statistic $T_k$ should behave in the same fashion when applied to *any* stretch of size $k$ of consecutive data points extracted from the data series $X_1, \ldots, X_n$; this observation motivates the notion of "scanning." On top of the particular application that will become obvious immediately, scanning may also provide an alternative way to think about the usual expanding sample asymptotics for stationary time series.

DEFINITION 3.1. A scan is a collection of $n$ block-subsamples of the sequence $X_1, \ldots, X_n$ with the following two properties: (a) within each scan there is a single block of each size $k = 1, \ldots, n$; and (b) those $n$ blocks are nested, that is, the block of size $k_1$ can be found as a stretch within the block of size $k_2$ when $k_1 \leq k_2$.

As usual, a block-subsample of the sequence $X_1, \ldots, X_n$ is a block of consecutive observations, that is, a set of the type $X_j, X_{j+1}, \ldots, X_{j+m}$.

We will say that the sequence $X_1, \ldots, X_n$ has been *scanned* if a block corresponding to each block size $k = 1, \ldots, n$ has been extracted, and if those blocks are nested, that is, the block of size $k_1$ can be found as a stretch within the block of size $k_2$ when $k_1 \leq k_2$. For example, in Section 3.1 the following "*direct*" scan was employed:

$$(X_1), (X_1, X_2), (X_1, X_2, X_3), \ldots, (X_1, \ldots, X_{n-1}), (X_1, \ldots, X_n),$$

over which the in-sample "evolution" of $T_n$ was investigated. Nevertheless, there are many possible scans; for example, consider the "*reverse*" scan

$$(X_n), (X_{n-1}, X_n), (X_{n-2}, X_{n-1}, X_n), \ldots, (X_2, \ldots, X_n), (X_1, \ldots, X_n).$$

In general, a scan will start at time-point $j$ (say) and then the blocks will proceed growing/expanding to the left and/or to the right—thus, the different perspective on asymptotics; for example, a valid scan is

$$(X_5), (X_4, X_5), (X_3, X_4, X_5), (X_3, X_4, X_5, X_6), \ldots, (X_1, \ldots, X_n);$$

note how within each block the natural time order is preserved, and how all scans end with the block containing the full data set. The number of possible scans is large as the following proposition shows.

PROPOSITION 3.1. *There are $2^{n-1}$ different scans of the sequence $X_1, \ldots, X_n$ when no ties are present.*

Let $B_i^k = (X_i, \ldots, X_{i+k-1})$, that is, $B_i^k$ for $i = 1, \ldots, n - k + 1$ are all the possible blocks of size $k$. Pascal's triangle and a backward induction argument suggest the following useful corollary.

COROLLARY 3.1. *Among the $2^{n-1}$ different scans of the sequence $X_1, \ldots, X_n$, there are exactly $\binom{n-k}{i-1} 2^{k-1}$ scans that contain block $B_i^k$ as their block of size $k$ for $1 \leq i \leq n - k + 1$.*

A collection of algorithms to generate randomly selected scans can be found at www.math.ucsd.edu/~politis/PAPER/scansAlgorithms.pdf, where some properties of those algorithms are also discussed.

3.3. *Improving upon the basic estimator.* As mentioned before, the usage of the particular "direct" scan

$$(X_1), (X_1, X_2), (X_1, X_2, X_3), \ldots, (X_1, \ldots, X_{n-1}), (X_1, \ldots, X_n)$$

in Section 3.1 was quite arbitrary; any scan could have been used with similar results. To elaborate, consider all the $2^{n-1}$ different scans of the sequence $X_1, \ldots, X_n$; order the scans in some arbitrary fashion, focus on the $I$th such scan, and consider the following analogs of our previous statements (a)–(c).

(a[$I$]) Let $T_n = T_n(X_1, \ldots, X_n)$ be some positive statistic whose rate of convergence/divergence depends on some unknown real-valued parameter $\lambda$.

(b[$I$]) For $k = 1, \ldots, n$, let $T_k^{(I)}$ denote the value of the statistic $T_k$ as computed from the block of size $k$ of the $I$th scan of the sequence $X_1, \ldots, X_n$.

(c[$I$]) Estimate $\lambda$ by $\hat{\lambda}^{(I)} = g^{-1}(\hat{g})$, or by $\breve{\lambda}^{(I)} = g^{-1}(\breve{g})$, where

$$\hat{g} = \frac{\sum_{k=1}^{n}(Y_k - \bar{Y})(\log k - \overline{\log n})}{\sum_{k=1}^{n}(\log k - \overline{\log n})^2}, \qquad \breve{g} = \frac{\sum_{k=1}^{n} Y_k \log k}{\sum_{k=1}^{n} \log^2 k},$$

and $Y_k = \log T_k^{(I)}$ for $k = 1, \ldots, n$, $\bar{Y} = \frac{1}{n} \sum_{k=1}^{n} Y_k$ and $\overline{\log n} = \frac{1}{n} \sum_{k=1}^{n} \log k$.

THEOREM 3.2. *Assume that the time series $\{X_t\}$ is strictly stationary. Under the assumptions of Theorem 3.1, the conclusions of Theorem 3.1 remain true with $\hat{\lambda}^{(I)}$ and $\breve{\lambda}^{(I)}$ in place of $\hat{\lambda}$ and $\breve{\lambda}$, respectively, for any $I$.*

Theorem 3.2—whose proof is identical to the proof of Theorem 3.1—suggests an approach on potentially improving the estimators $\hat{\lambda}$ and $\breve{\lambda}$ by combining/averaging the estimators based on scans. Consider the estimators $\hat{\lambda}^{(1)}, \ldots, \hat{\lambda}^{(N)}$ and $\breve{\lambda}^{(1)}, \ldots, \breve{\lambda}^{(N)}$ for some integer $N$, and define $\hat{\lambda}^* =$ $N^{-1} \sum_{i=1}^{N} \hat{\lambda}^{(i)}$ and $\breve{\lambda}^* = N^{-1} \sum_{i=1}^{N} \breve{\lambda}^{(i)}$. A different way of combining estimators is given by the median; so, we also define $\hat{\lambda}^{**} = \text{median}\{\hat{\lambda}^{(1)}, \ldots, \hat{\lambda}^{(N)}\}$ and $\breve{\lambda}^{**} = \text{median}\{\breve{\lambda}^{(1)}, \ldots, \breve{\lambda}^{(N)}\}$. The median estimators $\hat{\lambda}^{**}$ and $\breve{\lambda}^{**}$ will exhibit similar variance reduction behavior as the mean estimators $\hat{\lambda}^*$ and $\breve{\lambda}^*$. However, the median may be preferable in practice because of its robustness. The following corollary shows that averaging does not hurt asymptotically.

COROLLARY 3.2. *Assume that the time series $\{X_t\}$ is strictly stationary.*

(i) *Assume $N$ is fixed. Under the assumptions of Theorem 3.1, the conclusions of Theorem 3.1 remain true with $\hat{\lambda}^*$ or $\hat{\lambda}^{**}$ in place of $\hat{\lambda}$, and $\breve{\lambda}^*$ or $\breve{\lambda}^{**}$ in place of $\breve{\lambda}$.*

(ii) *Assume $N$ is a general positive function of $n$ (possibly diverging to infinity as $n \to \infty$). Under the assumptions of Theorem 3.1, the conclusions of Theorem 3.1 remain true with $\hat{\lambda}^*$ in place of $\hat{\lambda}$ and $\breve{\lambda}^*$ in place of $\breve{\lambda}$.*

It is generally difficult to quantify the variance reduction effect of scanning estimators; nevertheless, the simulations in Section 6 show a very spectacular effect even with a small value of $N$. Note that $N$ is really tied to the practitioner's computational facilities, and not so much to the sample size $n$ or the number of scans $2^{n-1}$. The recommendation is to take $N$ as big as computationally feasible; in practice, however, even taking $N$ as small as 100 gives a significant benefit especially if the $N$ scans under consideration are very different from one another. A way to ensure this is to use $N$ randomly selected scans from an algorithm that gives (close to) equal weight to each scan. A practical option is given by Algorithm $A(f)$ or Algorithm $B'$—the latter being valid only for weakly dependent, stationary sequences; see www.math.ucsd.edu/~politis/PAPER/scansAlgorithms.pdf for details.

## 4. Extensions of the basic methodology.

4.1. *Limit theorems with centering.* As mentioned in Remark 3.2, centering can typically be omitted in the case of diverging statistics. By contrast, in most cases of converging statistics a centering will be necessary in order to transform $T_n$ into a bounded random variable (in probability). Therefore, the following extension of the rate estimation methodology of Section 3 is proposed.

(a′) Let $T_n = T_n(X_1, \ldots, X_n)$ be some (not necessarily positive) statistic whose rate of convergence depends on some unknown real-valued parameter $\lambda$. Also assume that $P(T_k = T_n) = 0$ for $k = 1, \ldots, n-1$.

(b′) Assume that for some slowly varying function $L(n) > 0$ and for some known invertible function $g(\cdot)$ that is continuous over an interval that contains $\lambda$, and such that $g(\lambda) < 0$, we have

$$n^{-g(\lambda)} L(n) |T_n - \mu| \stackrel{\mathcal{L}}{\Longrightarrow} J \qquad \text{as } n \to \infty, \tag{8}$$

where $\mu$ is a real-valued parameter and $J$ some well-defined probability distribution; both $\mu$ and the shape of the limit distribution $J$ can be unknown.

(c′) Let $m, b$ be positive integers with $m \leq n - b$ and $b \leq n$; as before, we estimate $g(\lambda)$ by $\hat{g}_{m,b} = \sum_{k=m}^{b+m}(Y_k - \bar{Y})(\log k - \overline{\log}) / \sum_{k=m}^{b+m}(\log k - \overline{\log})^2$, and $\lambda$ by $\hat{\lambda}_{m,b} = g^{-1}(\hat{g}_{m,b})$, where $Y_k = \log |T_k - T_n|$, $\bar{Y} = \frac{1}{b+1}\sum_{k=m}^{b+m} Y_k$ and $\overline{\log} = \frac{1}{b+1}\sum_{k=m}^{b+m} \log k$.

Note that $\hat{g}$ in the above is an $L_2$ regression estimator of slope. As in Remark 3.1, here too it should be stressed that an $L_1$ estimator of slope in the regression of $Y_k$ on $\log k$ for $k = m, \ldots, b + m$ (with an intercept term included) might well give an attractive alternative that would be robust to the possibility that one of the $T_k$'s happens to be very close to $T_n$.

THEOREM 4.1. *If statements* (a′), (b′) *and* (c′) *are true, and assumptions* (4) *and* (6) *hold, then* $\hat{\lambda}_{m,b} \stackrel{P}{\longrightarrow} \lambda$, *provided* $1 \leq m \leq n - b$ *and* $b \to \infty$ *but* $b + m = o(n)$ *as* $n \to \infty$.

The assumption $P(T_k = T_n) = 0$ is imposed to ensure that $Y_k$ is well defined; it follows easily if the distribution of the statistic $T_n$ is absolutely continuous, in which case the probability of exact ties is zero. The condition $P(T_k = T_n) = 0$ could actually be relaxed to $P(T_k = T_n) \to 0$ when $k = k(n) \to \infty$ as $n \to \infty$ to accommodate the handling of statistics with discrete distributions; the details are straightforward and are omitted.

REMARK 4.1. Note that choosing $m$ is not a "bandwidth" selection problem; the choice $m = 1$ is fine for Theorem 4.1, although, in practice, one may prefer to take $m$ to be a small positive integer. Nevertheless, the trade-off requirements $b \to \infty$ but $b + m = o(n)$ imply that choosing $b$ is unfortunately a "bandwidth"-type problem. In this sense, rate estimation for uncentered diverging statistics seems to be easier to deal with; see, for example, Remark 3.2. To sidestep this difficulty, one may try to recast the problem into a diverging setup. So if $T_n$ is nonnegative, and if a lower bound for $g(\lambda)$ is known to exist [say $G < g(\lambda) < 0$], then line (8) implies that the uncentered quantity $n^{G-g(\lambda)} L(n) T_n$ should be diverging to $\infty$, and thus, the methods of Section 3 may be applicable; see Section 5.1 for an example of such a transformation.

4.2. *Improving upon the basic estimator.* As before, the notion of scanning may lead to improved estimation. Focus on the $I$th scan, and let $T_k^{(I)}$ denote the value of $T_k$ computed from the block of size $k$ of the $I$th scan of the sequence $X_1, \ldots, X_n$. Estimate $g(\lambda)$ by

$$\hat{g}_{m,b} = \frac{\sum_{k=m}^{b+m}(Y_k - \bar{Y})(\log k - \overline{\log})}{\sum_{k=m}^{b+m}(\log k - \overline{\log})^2},$$

and $\lambda$ by $\hat{\lambda}_{m,b}^{(I)} = g^{-1}(\hat{g}_{m,b})$, where $Y_k = \log|T_k^{(I)} - T_n|$, $\bar{Y} = \frac{1}{b+1}\sum_{k=m}^{b+m} Y_k$ and $\overline{\log} = \frac{1}{b+1}\sum_{k=m}^{b+m} \log k$. The following theorem and corollary ensue with proof identical to the proof of Theorem 4.1 combined with Corollary 3.2.

THEOREM 4.2. *Assume the time series $\{X_t\}$ is strictly stationary. If statements* (a$'$), (b$'$) *and* (c$'$) *are true, and assumptions* (4) *and* (6) *hold, then $\hat{\lambda}_{m,b}^{(I)} \xrightarrow{P} \lambda$, provided $1 \le m \le n - b$ and $b \to \infty$ but $b + m = o(n)$ as $n \to \infty$.*

To produce an improved estimator, we may again define

$$\hat{\lambda}_{m,b}^* = N^{-1}\sum_{i=1}^{N}\hat{\lambda}_{m,b}^{(i)} \quad \text{and} \quad \hat{\lambda}_{m,b}^{**} = \text{median}\{\hat{\lambda}_{m,b}^{(1)}, \ldots, \hat{\lambda}_{m,b}^{(N)}\},$$

where $N$ is some fixed positive integer.

COROLLARY 4.1. *Assume the time series $\{X_t\}$ is strictly stationary. If statements* (a$'$), (b$'$) *and* (c$'$) *are true, and assumptions* (4) *and* (6) *hold, then $\hat{\lambda}_{m,b}^* \xrightarrow{P} \lambda$ and $\hat{\lambda}_{m,b}^{**} \xrightarrow{P} \lambda$, provided $1 \le m \le n - b$ and $b \to \infty$ but $b + m = o(n)$ as $n \to \infty$.*

**5. Two examples with long memory.** The study of long memory time series appears to have been initiated by the hydrologist H. E. Hurst [12], who investigated the flow of the river Nile. Notably, Hurst's original R/S statistic was driven by a log-log regression as is our rate estimator $\hat{\lambda}$; see Beran [2] or Giraitis, Robinson and Surgailis [10] and the references therein. Interestingly, the well-known Geweke and Porter-Hudak [9] estimator of the long memory parameter also entails a log-regression based on some particular diverging statistics, namely, the periodogram ordinates at frequencies near zero.

5.1. *A second example*: *long memory time series.* Long memory time series are typically defined via an underlying stationary, mean zero, purely nondeterministic Gaussian time series $\{G_t, t \in \mathbf{Z}\}$ with autocovariance $R(k) = \text{Cov}(G_0, G_k)$ that is not absolutely summable. So assume that $X_t = h(G_t)$, where $h$ is some measurable function satisfying $Eh^2(G_t) < \infty$. Also assume that $R(k) = k^{-\beta}L(k)$ as $k \to \infty$, where $L$ is some slowly varying function, and $\beta > 0$ some unknown constant termed the long memory parameter. If $\beta > 1$, then the series $\{X_t\}$ and $\{G_t\}$ are weakly dependent, and the following central limit theorem typically holds:

$$\sqrt{n}(\bar{X}_n - \mu) \stackrel{\mathcal{L}}{\Longrightarrow} \mathrm{N}\left(0, \sum_{k=-\infty}^{\infty} R_X(k)\right), \tag{9}$$

where $\bar{X}_n = n^{-1}\sum_{t=1}^n X_t$, $R_X(k) = \mathrm{Cov}(X_0, X_k)$ and $\mu = EX_t$. If $\beta \leq 1$, then the sequences $\{X_t\}$ and $\{G_t\}$ are said to be *long-range dependent* and neither of them is strong mixing; see Ibragimov and Rozanov [13]. Hence, the subsampling methodology of Politis and Romano [18] may not be applicable, and the same is true for the subsampling rate estimator of Bertail, Politis and Romano [3] and Politis, Romano and Wolf [19], Chapter 8. In the long-range dependence case of $\beta \leq 1$, the following is true:

$$n(\bar{X}_n - \mu)/d_n \stackrel{\mathcal{L}}{\Longrightarrow} W_q, \tag{10}$$

where $d_n = n^{1-q\beta/2}L^{q/2}(n)$, and $q$ is the Hermite rank of $h$; see Taqqu [20, 21]. It is often the case that $q = 1$, in which case the limit distribution $W_1$ is a mean-zero Gaussian; for $q \geq 2$, $W_q$ is not Gaussian. Nevertheless, just the existence of the limit distributions in lines (9) and (10) is enough to imply that the techniques of Section 4 are applicable. In particular, a consistent estimator of the product $q\beta$ can be constructed using the sample mean as the converging statistic in Theorem 4.1; if $q$ is known, then this immediately yields an estimator of the long memory parameter $\beta$.

Different statistics could also be used; one example is the familiar second sample moment $S_n^2 = n^{-1}\sum_{t=1}^n X_t^2$ that was the focus of Section 2. As analogous limit theorems as (9) and (10) hold for the second sample moment, our rate estimation method of Theorem 4.1 could be based on the converging statistic $S_n^2$. The second sample moment $S_n^2$, however, is especially useful as it can be transformed to a diverging statistic as suggested in Remark 4.1. To do this, we simply let $T_n = \sum_{t=1}^n X_t^2 = nS_n^2$. It is easy to see that the requirements of Theorem 3.1 are satisfied for the diverging statistic $T_n$, and thus, a "bandwidth-free," consistent estimator of the product $q'\beta$ can be built based on $T_n$; here, of course, $q'$ denotes the Hermite rank of the function $h^2$.

5.2. *A third example*: *heavy tails with long memory*. Consider a time series defined as $X_t = \sqrt{\varepsilon_t}G_t$ for $t \in \mathbf{Z}$, where the series $\{\varepsilon_t\}$ and $\{G_t\}$ are independent, and the $\varepsilon_t$'s are positive and i.i.d. with distribution in $D(\alpha/2)$ for some $\alpha \in (0,2)$, and $\{G_t\}$ is stationary Gaussian with mean zero and autocovariance $R(k)$. For some $\zeta \in [0,1)$, define the condition

$$LM(\zeta): \left\{ \sum_{|h|<n} R(h) \sim Cn^\zeta \text{ and } \sum_{|h|<n} |R(h)| = O(n^\zeta) \text{ as } n \to \infty \right\},$$

where $C > 0$ is a constant. As before, the series $\{X_t\}$ and $\{Z_t\}$ are said to have long memory if $LM(\zeta)$ holds with $\zeta \in (0,1)$, in which case the long memory parameter $\beta$ equals $1-\zeta$; the case $LM(0)$ denotes weak dependence.

Interestingly, when appropriately normalized, the sample second moment converges in distribution in this general setting as the following proposition demonstrates; see Gomes, de Haan and Pestana [11] and McElroy and Politis [15] for related results.

PROPOSITION 5.1. *In the setting described above [including condition $LM(\zeta)$], suppose that $\varepsilon_t$ is absolutely continuous with a probability density $f_\varepsilon$ that is bounded and ultimately monotone, that is, $f_\varepsilon$ is monotone on $(z, \infty)$ for some $z > 0$, and is monotone on $(-\infty, u)$ for some $u < 0$. Then we have*

$$a_n^{-2} \sum_{t=1}^n X_t^2 \overset{\mathcal{L}}{\Longrightarrow} W \qquad \text{as } n \to \infty,$$

*where $a_n = n^{1/\alpha} K(n)$ for some slowly varying function $K(n)$. In the above, $W$ is $\alpha/2$-stable with scale $C_{\alpha/2}^{-2/\alpha}(E|G_t|^\alpha)^{2/\alpha}$, skewness 1 and location zero, and the constants $C_p^{-1}$ are defined by $C_p^{-1} = \Gamma(2-p)\cos(\pi p/2)/(1-p)$.*

The limit theorem of Proposition 5.1 is interesting, because the convergence of the sample second moment does *not* depend on the long memory parameter, and hence, our methods from Sections 3 and 4 can be unambiguously applied to estimate $\alpha$. Other methods in the tail index estimation literature may well encounter serious difficulties in this context, being sensitive to long-range dependence; this seems to be true for the Hill estimator—see Embrechts, Klüppelberg and Mikosch [8]. It is also true for the subsampling rate estimator of Bertail, Politis and Romano [3]; see the discussion in Section 5.1.

**6. A small simulation experiment.** We now revisit the setup of Section 2, that is, data $X_1, \ldots, X_n$ from a linear time series; see Remark 3.3. To perform the simulation, the AR(1) model

(11) $$X_t = \rho X_{t-1} + Z_t$$

was employed with $\rho = -0.5$, 0.1 or 0.7 and $\{Z_t\}$ i.i.d. from a distribution $F \in D(\alpha)$. The distributions used were (i) $\{Z_t\} \sim$ i.i.d. Cauchy, (ii) $\{Z_t\} \sim$

i.i.d. 1.5–Stable (symmetric), (iii) $\{Z_t\} \sim$ i.i.d. 1.9–Stable (symmetric), (iv) $\{Z_t\} \sim$ i.i.d. N(0, 1), (v) $\{Z_t\} \sim$ i.i.d. Pareto(2, 1), (vi) $\{Z_t\} \sim$ i.i.d. Burr (2, 1, 0.5) and (vii) $Z_t = \tilde{Z}_t \cdot \max(1, \log_{10}|\tilde{Z}_t|)$, where $\{\tilde{Z}_t\} \sim$ i.i.d. Burr (2, 1, 0.5). The variation (vii) has as its purpose the construction of a nonnormal domain of attraction, that is, the case where the slowly-varying function $L$ is not constant; see Embrechts, Klüppelberg and Mikosch [8].

For each combination of the value of $\rho$ and the distribution $F$, 100 time series stretches were generated, each of length $n = 1{,}000$. From each series, the estimator $\hat{\alpha}$ was computed, where $\hat{\alpha}$ was defined in Section 2; also computed were the improved versions $\hat{\alpha}^*$ and $\hat{\alpha}^{**}$, that is, the mean and median of the values of $\hat{\alpha}$ based on scans as in Corollary 3.2. Note that the information that $0 < \alpha \le 2$ was explicitly used in that values of $\hat{\alpha}$ bigger than 2 were truncated to the value 2; interestingly, no occurrences of a negative $\hat{\alpha}$ were observed. This truncation is necessary for good performance of $\hat{\alpha}^*$, but is superfluous for $\hat{\alpha}^{**}$ since the latter is based on a median that "clips" outliers.

A number of scanning algorithms can be devised; the website www.math.ucsd.edu/~politis/PAPER/scansAlgorithms.pdf presents Algorithms A, B and A($f$), making the claim that Algorithm A($f$)—with a carefully chosen $f$—may be preferable. However, Algorithm A($f$) is very computer-intensive. Although this is not a problem for the practitioner with a single dataset at hand, it is prohibitive in terms of conducting a simulation with thousands of datasets. A computational shortcut is presented by Algorithm B' that is valid for weakly dependent stationary sequences only. In particular, it is *not* suitable for the "long-memory" series of Section 5; see the aforementioned website for more details.

The results of our simulation, where $N$ random scans were generated using Algorithm B', are summarized in Table 1 where the empirical mean squared error (MSE) of each estimator is given. In this setup, the benchmark for comparison among estimators of $\alpha$ is given by the Hill estimator $H_q$ based on $q$ extreme order statistics. Empirical MSEs of Hill estimators are given in Table 2 for different values of $q$. Also included in Table 2 are the (empirically found) true optimal values of $q$, denoted by $q_{\text{opt}}$; in other words, $H_{q_{\text{opt}}}$ was the smallest MSE empirically computed from the model in question over a wide range of $q$ values. Things to note are the following:

- Averaging over scans does indeed succeed in dramatically reducing the MSE of estimation. As a matter of fact, even with $N$ as low as 100, significant benefits ensue, typically halving the MSE of the original estimator; this is of course contingent on having those $N$ scans generated in a very "random" fashion as Algorithm B' ensures.
- The comparison between $\hat{\alpha}^*$ and $\hat{\alpha}^{**}$ is unclear. The former seems to lead to somewhat smaller MSEs, but it should be borne in mind that its performance is aided by the truncation of the original estimator to the value 2. On the other hand, $\hat{\alpha}^{**}$ is more robust, and thus recommendable in a general setup when outside information—such as the restriction $\alpha \in (0, 2]$—may not be available.

TABLE 1
*Empirical MSEs of estimators of the heavy tail index $\alpha$; data from model (11) with $n = 1{,}000$ and (a) $\rho = 0.1$, (b) $\rho = 0.7$, (c) $\rho = -0.5$*

|  | $\hat{\alpha}$ | $\hat{\alpha}^*_{(N=20)}$ | $\hat{\alpha}^*_{(N=100)}$ | $\hat{\alpha}^*_{(N=200)}$ | $\hat{\alpha}^{**}_{(N=20)}$ | $\hat{\alpha}^{**}_{(N=100)}$ | $\hat{\alpha}^{**}_{(N=200)}$ |
|---|---|---|---|---|---|---|---|
| (a) | | | | | | | |
| (i) | 0.315 | 0.223 | 0.102 | 0.096 | 0.329 | 0.098 | 0.085 |
| (ii) | 0.171 | 0.109 | 0.064 | 0.064 | 0.152 | 0.107 | 0.109 |
| (iii) | 0.051 | 0.036 | 0.025 | 0.024 | 0.044 | 0.041 | 0.037 |
| (iv) | 0.006 | 0.004 | 0.002 | 0.002 | 0.004 | ($<0.0005$) | ($<0.0005$) |
| (v) | 0.222 | 0.190 | 0.142 | 0.140 | 0.220 | 0.167 | 0.166 |
| (vi) | 0.294 | 0.159 | 0.079 | 0.079 | 0.228 | 0.106 | 0.101 |
| (vii) | 0.319 | 0.156 | 0.074 | 0.068 | 0.260 | 0.106 | 0.096 |
| (b) | | | | | | | |
| (i) | 0.328 | 0.193 | 0.108 | 0.106 | 0.265 | 0.127 | 0.109 |
| (ii) | 0.161 | 0.097 | 0.057 | 0.055 | 0.147 | 0.101 | 0.093 |
| (iii) | 0.078 | 0.046 | 0.034 | 0.033 | 0.059 | 0.058 | 0.052 |
| (iv) | 0.011 | 0.009 | 0.006 | 0.005 | 0.010 | 0.002 | 0.001 |
| (v) | 0.120 | 0.079 | 0.079 | 0.077 | 0.091 | 0.088 | 0.084 |
| (vi) | 0.343 | 0.205 | 0.105 | 0.103 | 0.312 | 0.112 | 0.107 |
| (vii) | 0.314 | 0.189 | 0.062 | 0.060 | 0.295 | 0.102 | 0.097 |
| (c) | | | | | | | |
| (i) | 0.322 | 0.234 | 0.145 | 0.138 | 0.322 | 0.156 | 0.145 |
| (ii) | 0.139 | 0.097 | 0.055 | 0.052 | 0.151 | 0.091 | 0.086 |
| (iii) | 0.054 | 0.046 | 0.026 | 0.028 | 0.056 | 0.040 | 0.044 |
| (iv) | 0.008 | 0.005 | 0.003 | 0.003 | 0.007 | 0.001 | ($<0.0005$) |
| (v) | 0.254 | 0.193 | 0.164 | 0.169 | 0.218 | 0.202 | 0.210 |
| (vi) | 0.295 | 0.204 | 0.089 | 0.079 | 0.283 | 0.123 | 0.109 |
| (vii) | 0.321 | 0.151 | 0.064 | 0.056 | 0.237 | 0.105 | 0.097 |

Comparing Table 1 to Table 2, it is apparent that both $\hat{\alpha}^*$ and $\hat{\alpha}^{**}$ underperform as compared to the optimized Hill estimator $H_{q_{\text{opt}}}$ in cases (i), (ii), (v) and (vi), whereas $\hat{\alpha}^*$ and $\hat{\alpha}^{**}$ perform comparably to $H_{q_{\text{opt}}}$ in cases (iii) and (vii). Both $\hat{\alpha}^*$ and $\hat{\alpha}^{**}$ perform excellently in the Gaussian case (iv); however, as kindly pointed out by one of the referees, the Hill estimator is inapplicable/inconsistent in this case as it diverges to infinity—therefore, the n/a's in Table 2.

Perhaps it should be stressed that $q_{\text{opt}}$ is *not* known by the practitioner. As mentioned earlier, estimation of $q_{\text{opt}}$ is not a trivial matter and is further complicated when the data are dependent; see Embrechts, Klüppelberg and

TABLE 2
*Empirical MSEs of Hill estimator $H_q$ based on $q$ order statistics; data from model (11) with $n = 1{,}000$ and (a) $\rho = 0.1$, (b) $\rho = 0.7$, (c) $\rho = -0.5$*

|       | $H_{100}$ | $H_{200}$ | $H_{300}$ | $H_{400}$ | $H_{q_{\mathrm{opt}}}$ | $q_{\mathrm{opt}}$ |
|-------|-----------|-----------|-----------|-----------|------------------------|--------------------|
| (a)   |           |           |           |           |                        |                    |
| (i)   | 0.011     | 0.011     | 0.048     | 0.170     | 0.007                  | 140                |
| (ii)  | 0.121     | 0.013     | 0.130     | 0.546     | 0.013                  | 200                |
| (iii) | 1.469     | 0.043     | 0.290     | 1.147     | 0.017                  | 220                |
| (iv)  | n/a       | n/a       | n/a       | n/a       | n/a                    | n/a                |
| (v)   | 0.149     | 0.291     | 0.450     | 0.629     | 0.086                  | 40                 |
| (vi)  | 0.032     | 0.065     | 0.099     | 0.138     | 0.027                  | 60                 |
| (vii) | 0.094     | 0.106     | 0.134     | 0.167     | 0.059                  | 20                 |
| (b)   |           |           |           |           |                        |                    |
| (i)   | 0.045     | 0.019     | 0.051     | 0.136     | 0.019                  | 200                |
| (ii)  | 0.253     | 0.039     | 0.135     | 0.562     | 0.031                  | 220                |
| (iii) | 1.262     | 0.050     | 0.315     | 1.198     | 0.035                  | 220                |
| (iv)  | n/a       | n/a       | n/a       | n/a       | n/a                    | n/a                |
| (v)   | 0.373     | 0.147     | 0.059     | 0.026     | 0.026                  | 400                |
| (vi)  | 0.057     | 0.023     | 0.017     | 0.022     | 0.017                  | 300                |
| (vii) | 0.048     | 0.064     | 0.078     | 0.087     | 0.048                  | 100                |
| (c)   |           |           |           |           |                        |                    |
| (i)   | 0.017     | 0.012     | 0.042     | 0.155     | 0.011                  | 180                |
| (ii)  | 0.135     | 0.015     | 0.118     | 0.532     | 0.013                  | 220                |
| (iii) | 1.297     | 0.042     | 0.286     | 1.111     | 0.015                  | 220                |
| (iv)  | n/a       | n/a       | n/a       | n/a       | n/a                    | n/a                |
| (v)   | 0.184     | 0.421     | 0.727     | 1.072     | 0.118                  | 40                 |
| (vi)  | 0.038     | 0.084     | 0.141     | 0.219     | 0.034                  | 60                 |
| (vii) | 0.104     | 0.138     | 0.183     | 0.252     | 0.052                  | 20                 |

Mikosch [8] or Danielsson et al. [7] and the references therein. This phenomenon is manifested in our simulations, especially in cases (v)–(vii), that is, the Pareto and Burr distributions, for which the value of the empirically found $q_{\mathrm{opt}}$ seems to be quite unstable as a function of the dependence factor $\rho$, to the extent exemplified in our small simulation.

The simulation confirms that our proposed methodology leads to reasonable estimates of the index of domain of attraction under (linear) dependence and possibly nonnormal domains of attraction, that is, nonconstant slowly-varying function $L$. Nevertheless, it should be stressed that our methodology has general applicability, and it is not specific to the particular context as Hill's estimator is. Of course, it is expected that context-specific, carefully optimized estimators may give improved performance relative to this general "off-the-shelf" tool. The fact that in some of the cases considered, for example, (iii) and (vii), our general methodology performs comparably with the optimally fine-tuned Hill estimator (using the true $q_{\text{opt}}$) can be considered remarkable.

An added bonus of our methodology is that it is totally automatic: no fine-tuning is required in terms of a tricky "bandwidth"-type choice, such as estimating $q_{\text{opt}}$ for the Hill estimator. In addition, note that—even in the specific tail estimation context of this section—our methodology is applicable in connection with different diverging statistics other than the second moment. There is a plethora of such diverging estimators that can be used; for example, $T_n$ could be taken as the $2r$th sample moment for some integer $r \geq 1$, the $r$th sample moment of the absolute values of the $X_t$'s for some integer $r \geq 2$, the maximum $M_n = \max\{X_1, \ldots, X_n\}$ or the range $K_n = M_n - L_n$, where $L_n = \min\{X_1, \ldots, X_n\}$.

The performance of those different candidate statistics is context-specific, and will generally depend on many factors, including the underlying value of $\alpha$ as well. Furthermore, since all these different statistics yield useful information for $\alpha$, it is conceivable that they can all be combined to construct an improved estimator. To give a concrete example, let $\hat{\alpha}^{**(r)}$ denote our median-averaged estimator of $\alpha$ based on the $r$th sample moment of the absolute values as the diverging statistic. The estimators $\hat{\alpha}^{**(r)}$ for $r = 2, 3, \ldots, R$ can be constructed for some fixed integer $R$ whose magnitude will depend on the practitioner's computational facilities. Those $R$ estimators can then be combined to yield the yet improved estimator

$$\hat{\alpha}^{**,R} = \text{median}(\hat{\alpha}^{**(2)}, \ldots, \hat{\alpha}^{**(R)}). \tag{12}$$

## APPENDIX: TECHNICAL PROOFS

PROOF OF PROPOSITION 2.1. By a corollary of Karamata's representation theorem, see, for example, Theorem A.3.3 in Embrechts, Klüppelberg and Mikosch [8], it follows that $\log L(n)/\log n$ asymptotically behaves as $\int_z^n (\delta(u)/u)\, du/\log n$ for some number $z > 0$ and a measurable function $\delta(u)$ that tends to zero as $u \to \infty$. If the integral converges, then the assertion is proved; otherwise use l'Hôpital's rule to obtain an asymptotic rate of $\delta(n)$, which tends to zero. Thus, $\log L(n)/\log n = o(1)$. □

PROOF OF THEOREM 3.1. Consider the identity

$$Y_k = g(\lambda) \log k + U_k - \log L(k) \tag{13}$$

for $k = 1, \ldots, n$. After some straightforward calculations, we have that $\breve{g} = g(\lambda) + a_1 + a_2$, where

$$a_1 = \frac{\sum_{k=1}^n U_k \log k}{\sum_{k=1}^n \log^2 k} \quad \text{and} \quad a_2 = -\frac{\sum_{k=1}^n \log L(k) \log k}{\sum_{k=1}^n \log^2 k}.$$

Note that

$$c_1 \log n \leq \overline{\log n} \leq c_2 \log n \quad \text{and} \quad c_3 (\log n)^2 \leq n^{-1} \sum_{k=1}^{n} (\log k)^2 \leq c_4 (\log n)^2,$$

for some constants $c_i > 0$. Since it is assumed that $U_n = O_P(1)$, it follows that $a_1 = O_P(1/\log n) = o_P(1)$. Using line (2), it follows that $a_2 = o(1)$ as well. Hence, $\breve{g} = g(\lambda) + o_P(1)$. Finally, part (i) is proven by an application of the continuous mapping theorem.

To analyze $\hat{\lambda}$, note that similarly we have $\hat{g} = g(\lambda) + A_0 + A_2$, where

$$A_0 = \frac{\sum_{k=1}^{n}(U_k - \bar{U})(\log k - \overline{\log n})}{\sum_{k=1}^{n}(\log k - \overline{\log n})^2},$$

$$A_2 = -\frac{\sum_{k=1}^{n}(\log L(k) - \overline{\log L})(\log k - \overline{\log n})}{\sum_{k=1}^{n}(\log k - \overline{\log n})^2};$$

here $\bar{U} = \frac{1}{n} \sum_{k=1}^{n} U_k$ and $\overline{\log L} = \frac{1}{n} \sum_{k=1}^{n} \log L(k)$. Let $A_1 = E[A_0]$. By a Riemann-sum approximation argument, it follows that

(14)
$$n^{-1} \sum_{k=1}^{n} (\log k - \overline{\log n})^2 = n^{-1} \sum_{k=1}^{n} (\log k/n - \overline{\log k/n})^2$$
$$\to \int_0^1 (\log x + 1)^2 \, dx = 1,$$

where $\overline{\log k/n} = \frac{1}{n} \sum_{k=1}^{n} \log k/n$. Focus on the numerator of $A_2$: defining $L_n(x) = L(\lceil nx \rceil)$ such that $L_n(k/n) = L(k)$, we obtain

$$-\frac{1}{n} \sum_{k=1}^{n} (\log L(k) - \overline{\log L})(\log k - \overline{\log n})$$

$$= -\frac{1}{n} \sum_{k=1}^{n} \log L_n(k/n)(\log k/n - \overline{\log k/n})$$

$$= -\int_0^1 L_n(x)(\log x + 1) \, dx + o(1)$$

by a straightforward application of the definition of the Riemann integral. Note that the error in this approximation is just $o(1)$ instead of $O(1/n)$, since $\log x$ does not have a bounded derivative on $[0,1]$. From Theorem A3.3 of Embrechts, Klüppelberg and Mikosh [8], we have the representation $L(y) = c(y) \exp\{\int_z^y (\eta(u)/u) \, du\}$ for some $z > 0$, $c(x) \to c > 0$ and $\eta(x) \to 0$ as $x \to \infty$. It follows that

$$\log L_n(x) - \log L_n(1) = \log(c_{\lceil nx \rceil}/c_n) - \int_{\lceil nx \rceil}^{n} (\eta(u)/u) \, du.$$

So for each fixed $x \in (0,1]$, $c_{\lceil nx \rceil}/c_n \to 1$ and

$$\left| \int_{\lceil nx \rceil}^{n} (\eta(u)/u) \, du \right| \leq \sup_{u \in [nx,n]} |\eta(u)| \frac{n(1-x)}{nx},$$

which tends to zero as $n \to \infty$. This shows that $\log L_n(x) - \log L_n(1) = o(1)$. But since $\log L_n(1)$ does not depend on $x$ and $\int_0^1 (\log x + 1) \, dx = 0$,

$$-\int_0^1 L_n(x)(\log x + 1) \, dx = -\int_0^1 (L_n(x) - L_n(1))(\log x + 1) \, dx \to 0$$

by the dominated convergence theorem, since the integrand converges uniformly to zero. This shows that $A_2 = o(1)$.

*Part* (ii). Now assume (4). From line (14), we also have that

$$A_0 = n^{-1} \sum_{k=1}^{n} U_k (\log k/n - \overline{\log k/n}) + o(1),$$

which we will denote by $I_1$. Now since $EU_n^2 \to EU^2$, we find that $\sup_n EU_n^2 < \infty$ so that $\{U_n\}$ is a uniformly integrable sequence. Together with $U_n \overset{\mathcal{L}}{\Longrightarrow} U$, this implies that $EU_n \to EU$ as $n \to \infty$, and also that for each $x$, $EU_n(x) \to EU$ with $U_n(x) = U_{\lceil nx \rceil}$. Hence we calculate

$$EI_1 = \frac{1}{n} \sum_{k=1}^{n} EU_n(k/n)(\log k/n - \overline{\log k/n})$$

$$= o(1) + \int_0^1 EU_n(x)(\log x + 1) \, dx$$

$$\to \int_0^1 EU(\log x + 1) \, dx = 0$$

using the dominated convergence theorem. Hence, $E\hat{g} = g(\lambda) + A_1 + A_2 = g(\lambda) + o(1)$. Now observe that

$$\text{Var}(A_0) \sim n^{-2} \sum_{k=1}^{n} \sum_{b=1}^{n} \text{Cov}(U_b, U_k)(\log b - \overline{\log n})(\log k - \overline{\log n}).$$

So from (4) it follows that $\text{Cov}(U_b, U_k) = O(1)$; thus, line (14) implies that $\text{Var}(A_0) = O(1)$, completing the proof of part (ii).

*Part* (iii). Now assume (5) as well;

$$EI_1 = \frac{1}{n} \sum_{k=1}^{n} (EU_k - EU)(\log k - \overline{\log n}),$$

since $EU$ does not depend on $k$. Taking absolute values produces a crude bound of $\frac{1}{n} \sum_{k=1}^{n} C k^{-p} \log k$ for some constant $C > 0$; thus, $A_1$ is clearly $O(n^{-p} \log n)$, and $A_2$ has already been analyzed.

*Part* (iv). Finally, assume (4) and (6). Then we have

$$\mathrm{Var}(A_0) \leq \frac{\log^2 n}{n^2} \sum_{k=1}^{n} \sum_{b=1}^{n} |\mathrm{Cov}(U_b, U_k)| = O(n^{\gamma_1 - \gamma_2}(\log^2 n)\tilde{L}(n)) = o(1),$$

which shows that $\hat{g} \xrightarrow{P} g(\lambda)$ and hence $\hat{\lambda} \xrightarrow{P} \lambda$ as well. $\square$

PROOF OF PROPOSITION 3.1. To see this, one has to look at the last block of a scan and deconstruct it, that is, go backward. Since the last block is always $(X_1, \ldots, X_n)$, the next-to-last is either $(X_1, \ldots, X_{n-1})$ or $(X_2, \ldots, X_n)$, that is, one of two choices. Similarly, from any step of the process, for example, from the block of size $k+1$, there are two choices for the preceding block corresponding to shrinking from the left or from the right. Thus, there are two choices for each of the $n-1$ steps of the deconstruction of the last blocks; these choices multiply to give the number $2^{n-1}$. $\square$

PROOF OF COROLLARY 3.2. Regarding $\hat{\lambda}^*$ and $\check{\lambda}^*$, the proof follows by a simple application of the Cauchy–Schwarz inequality. Regarding $\hat{\lambda}^{**}$ and $\check{\lambda}^{**}$, just note that they represent medians of $N$ i.i.d. random variables where $N$ is finite. $\square$

PROOF OF THEOREM 4.1. We first show that $\hat{\lambda}_{m,b} \xrightarrow{P} \lambda$ under the assumptions of Theorem 4.1 together with the additional assumption that $m \to \infty$. First note that if we define

$$U_k = \log(k^{-g(\lambda)} L(k)|T_k - T_n|) \qquad \text{for } k = m, \ldots, b + m,$$

then the identity (13) still holds true but now for $k = m, \ldots, b + m$ only. In addition, we also have $U_n = O_P(1)$ as $n \to \infty$ as before. To see this, note that $k^{-g(\lambda)} L(k)|T_k - \mu| \xrightarrow{\mathcal{L}} J$ as $k \to \infty$ by assumption (8). Also note

$$k^{-g(\lambda)} L(k)(T_k - \mu) = k^{-g(\lambda)} L(k)(T_k - T_n + T_n - \mu) = k^{-g(\lambda)} L(k)(T_k - T_n) + A_0,$$

where $A_0 = k^{-g(\lambda)} L(k)(T_n - \mu)$. But

$$|A_0| = k^{-g(\lambda)} L(k)|T_n - \mu| = \frac{k^{-g(\lambda)} L(k)}{n^{-g(\lambda)} L(n)} n^{-g(\lambda)} L(n)|T_n - \mu| = O_P\left(\left(\frac{k}{n}\right)^{-g(\lambda)}\right),$$

again by assumption (8). Since $k/n \leq (b+m)/n = o(1)$ and $g(\lambda) < 0$, it follows that $A_0 = o_P(1)$. Finally, Slutsky's theorem and the continuous mapping theorem ensure that $k^{-g(\lambda)} L(k)|T_k - T_n| \xrightarrow{\mathcal{L}} J$ as $k \to \infty$ and hence, $U_k = O_P(1)$ as $k \to \infty$. By a calculation similar to that in the proof of Theorem 3.1, we have that $\hat{g}_{m,b} = g(\lambda) + A_1 + A_2$, where now

$$A_1 = \frac{\sum_{k=m}^{b+m} (U_k - \bar{U})(\log k - \overline{\log})}{\sum_{k=m}^{b+m} (\log k - \overline{\log})^2}$$

and

$$A_2 = -\frac{\sum_{k=m}^{b+m}(\log L(k) - \overline{\log L})(\log k - \overline{\log})}{\sum_{k=m}^{b+m}(\log k - \overline{\log})^2};$$

here $\bar{U} = \frac{1}{b+1}\sum_{k=m}^{b+m} U_k$ and $\overline{\log L} = \frac{1}{b+1}\sum_{k=m}^{b+m} \log L(k)$. As in the proof of Theorem 3.1, it follows that, as $b \to \infty$, $A_2 = o(1)$ and $EA_1 = o(1)$ by equation (4). Moreover, $\operatorname{Var} A_1 = o(1)$ by equation (6), and hence, $\hat{g}_{m,b} = g(\lambda) + o_P(1)$. An application of the continuous mapping theorem shows that $\hat{\lambda}_{m,b} \xrightarrow{P} \lambda$.

We now wish to relax the extra assumption $m \to \infty$. To do this, we will show that $m=1$ is a good enough choice, that is, that $\hat{\lambda}_{1,b} \xrightarrow{P} \lambda$ when $b \to \infty$ but $b = o(n)$; the proof for other nondiverging choices for $m$ is similar. Note that by the above arguments we can write $\hat{\lambda}_{1,b} = g^{-1}(\hat{g}_{1,b})$, where

$$\hat{g}_{1,b} = g(\lambda) + A_1^* + A_2^*;$$

in the above, $A_1^*, A_2^*$ are similar to the terms $A_1, A_2$ but with summations of the type $\sum_{k=1}^{b+1}$ instead of $\sum_{k=m}^{b+m}$ in both numerator and denominator. Now consider a choice of $m$ satisfying $m \to \infty$ but also $m = o(b)$. By the above discussion, we have shown that $\hat{\lambda}_{m,b-m+1} \xrightarrow{P} \lambda$. In particular, we can write

$$\hat{g}_{m,b-m+1} = g(\lambda) + A_1' + A_2',$$

where $A_1', A_2'$ are again similar to the terms $A_1, A_2$ but with summations of the type $\sum_{k=m}^{b+1}$ instead of $\sum_{k=m}^{b+m}$; furthermore, we have also shown that $A_1', A_2'$ are both $o_P(1)$. Looking at the numerator of $A_1'$, we see a sum of the type $\sum_{k=m}^{b+1}$ which we have shown to be of order $O_P((b-m)\log(b-m))$. The denominator of $A_1'$ includes a sum of the type $\sum_{k=m}^{b+1}$ which is of exact order $O((b-m)\log^2(b-m))$. Now writing those sums as $\sum_{k=m}^{b+1} = \sum_{k=1}^{b+1} - \sum_{k=1}^{m-1}$ in both numerator and denominator of $A_1'$, and using the assumption $m = o(b)$, it follows that $A_1' = A_1^* + o_P(1)$. Similarly, $A_2' = A_2^* + o_P(1)$. Since $A_1', A_2'$ are both $o_P(1)$, it is immediate that $A_1^*, A_2^*$ are both $o_P(1)$; thus, $\hat{g}_{1,b} = g(\lambda) + o_P(1)$, and $\hat{\lambda}_{1,b} \xrightarrow{P} \lambda$ as desired. □

PROOF OF PROPOSITION 5.1. The proof is available at the website www.math.ucsd.edu/~politis/PAPER/scansAppendix1.pdf. □

**Acknowledgments.** Many thanks are due to Patrice Bertail (CREST and Univ. Paris X), Arup Bose (Indian Statistical Institute), Holger Drees (Univ. of Hamburg), Keith Knight (Univ. of Toronto), S. N. Lahiri (Iowa State Univ.), Ivan Mizera (Univ. of Alberta) and Joe Romano (Stanford) for helpful discussions, and to the Associate Editor and three anonymous referees for their constructive comments.

## REFERENCES

[1] Barbe, P. and Bertail, P. (1993). Testing the global stability of a linear model. Document de travail of INRA-CORELA, Paris. (See also its revised version, Document de travail CREST, Paris, 2003.)

[2] Beran, J. (1994). *Statistics for Long-Memory Processes*. Chapman and Hall, New York. MR1304490

[3] Bertail, P., Politis, D. N. and Romano, J. P. (1999). On subsampling estimators with unknown rate of convergence. *J. Amer. Statist. Assoc.* **94** 569–579. MR1702326

[4] Breiman, L. (1996). Bagging predictors. *Machine Learning* **24** 123–140.

[5] Csörgő, S., Deheuvels, P. and Mason, D. M. (1985). Kernel estimates of the tail index of a distribution. *Ann. Statist.* **13** 1050–1077. MR0803758

[6] Csörgő, S. and Viharos, L. (1998). Estimating the tail index. In *Asymptotic Methods in Probability and Statistics* (B. Szyszkowicz, ed.) 833–881. North-Holland, Amsterdam. MR1661521

[7] Danielsson, J., de Haan, L., Peng, L. and de Vries, C. G. (2001). Using a bootstrap method to choose the sample fraction in tail index estimation. *J. Multivariate Anal.* **76** 226–248. MR1821820

[8] Embrechts, P., Klüppelberg, C. and Mikosch, T. (1997). *Modelling Extremal Events for Insurance and Finance*. Springer, Berlin. MR1458613

[9] Geweke, J. and Porter-Hudak, S. (1983). The estimation and application of long memory time series models. *J. Time Ser. Anal.* **4** 221–238. MR0738585

[10] Giraitis, L., Robinson, P. M. and Surgailis, D. (1999). Variance-type estimation of long memory. *Stochastic Process. Appl.* **80** 1–24. MR1670119

[11] Gomes, M. I., de Haan, L. and Pestana, D. (2004). Joint exceedances of the ARCH process. *J. Appl. Probab.* **41** 919–926. MR2074832

[12] Hurst, H. E. (1951). Long-term storage capacity of reservoirs. *Trans. Amer. Soc. Civil Engineers* **116** 770–808.

[13] Ibragimov, I. A. and Rozanov, Y. A. (1978). *Gaussian Random Processes*. Springer, New York. MR0543837

[14] McElroy, T. and Politis, D. N. (2002). Robust inference for the mean in the presence of serial correlation and heavy-tailed distributions. *Econometric Theory* **18** 1019–1039. MR1926012

[15] McElroy, T. and Politis, D. N. (2007). Self-normalization for heavy-tailed time series with long memory. *Statist. Sinica* **17** 199–220.

[16] Meerschaert, M. and Scheffler, H. (1998). A simple robust estimation method for the thickness of heavy tails. *J. Statist. Plann. Inference* **71** 19–34. MR1651847

[17] Politis, D. N. (2002). A new approach on estimation of the tail index. *C. R. Math. Acad. Sci. Paris* **335** 279–282. MR1933674

[18] Politis, D. N. and Romano, J. P. (1994). Large sample confidence regions based on subsamples under minimal assumptions. *Ann. Statist.* **22** 2031–2050. MR1329181

[19] Politis, D. N., Romano, J. P. and Wolf, M. (1999). *Subsampling*. Springer, New York. MR1707286

[20] Taqqu, M. S. (1975). Weak convergence to fractional Brownian motion and to the Rosenblatt process. *Z. Wahrsch. Verw. Gebiete* **31** 287–302. MR0400329

[21] Taqqu, M. S. (1979). Convergence of integrated processes of arbitrary Hermite rank. *Z. Wahrsch. Verw. Gebiete* **50** 53–83. MR0550123

U.S. Bureau of the Census  
Statistical Research Division  
4700 Silver Hill Road  
Washington, DC 20233-9100  
USA  
E-mail: tucker.s.mcelroy@census.gov

Department of Mathematics  
University of California, San Diego  
La Jolla, California 92093-0112  
USA  
E-mail: dpolitis@ucsd.edu